# Confidence distribution (CD) – distribution estimator of a parameter

Kesar Singh[1,*], Minge Xie[1,†] and William E. Strawderman [1,‡]

*Rutgers University*

**Abstract:** The notion of confidence distribution (CD), an entirely frequentist concept, is in essence a Neymanian interpretation of Fisher's Fiducial distribution. It contains information related to every kind of frequentist inference. In this article, a CD is viewed as a distribution estimator of a parameter. This leads naturally to consideration of the information contained in CD, comparison of CDs and optimal CDs, and connection of the CD concept to the (profile) likelihood function. A formal development of a multiparameter CD is also presented.

## 1. Introduction and the concept

We are happy to dedicate this article to the memory of our colleague Yehuda Vardi. He was supportive of our efforts to develop this research area and in particular brought his paper with Colin Mallows (Mallows and Vardi [19]) to our attention during the discussion. A confidence-distribution (CD) is a compact expression of frequentist inference which contains information on just about every kind of inferential problem. The concept of a CD has its roots in Fisher's fiducial distribution, although it is a purely frequentist concept with a purely frequentist interpretation. Simply speaking, a CD of a univariate parameter is a data-dependent distribution whose $s$-th quantile is the upper end of a $100s\%$-level one-sided confidence interval of $\theta$. This assertion clearly entails that, for any $0 < s < t < 1$, the interval formed by $s$-th and $t$-th quantiles of a CD is a $100(t-s)\%$ level two-sided confidence interval. Thus, a CD is in fact Neymanian interpretation of Fisher's fiducial distribution (Neyman [21]). The concept of CD has appeared in a number of research articles. However, the modern statistical community has largely ignored the notion, particularly in applications. We suspect two probable causes lie behind this: (I) The first is its historic connection to Fisher's fiducial distribution, which is largely considered as "Fisher's biggest blunder" (see, for instance, Efron [8]); (II) Statisticians have not seriously looked at the possible utility of CDs in the context of modern statistical practice. As pointed out by Schweder and Hjort [22], there has recently been a renewed interest in this topic. Some recent articles include Efron [7, 8], Fraser [11, 12], Lehmann [15], Schweder and Hjort [22, 23], Singh, Xie and Strawderman [25, 26], among others. In particular, recent articles emphasize the Neymanian interpretation of the CD and present it as a valuable statistical tool for inference.

*Research partially supported by NSF DMS-0505552.
†Research partially supported by NSF SES-0241859.
‡Research partially supported by NSA Grant 03G-112.
[1]Department of Statistics, Hill Center, Busch Campus, Rutgers University, Piscataway, NJ 08854, USA, e-mail: kesar@stat.rutgers.edu; mxie@stat.rutgers.edu; starw@stat.rutgers.edu; url: stat.rutgers.edu

*AMS 2000 subject classifications:* primary 62F03, 62F12, 62G05; secondary 62G10, 62G20.

*Keywords and phrases:* confidence distribution, frequentist inference, fiducial inference, optimality, bootstrap, data-depth, likelihood function, p-value function.

For example, Schweder and Hjort [22] proposed reduced likelihood function from the CDs for inference, and Singh, Xie and Strawderman [26] developed attractive comprehensive approaches through the CDs to combining information from independent sources. The following quotation from Efron [8] on Fisher's contribution of the Fiducial distribution seems quite relevant in the context of CDs: "... but here is a safe prediction for the 21st century: statisticians will be asked to solve bigger and more complicated problems. I believe there is a good chance that objective Bayes methods will be developed for such problems, and that something like fiducial inference will play an important role in this development. Maybe Fisher's biggest blunder will become a big hit in the 21st century!"

In the remainder of this section, we give a formal definition of a confidence distribution and the associated notion of an asymptotic confidence distribution (aCD), provide a simple method of constructing CDs as well as several examples of CDs and aCDs. In the following formal definition of the CDs, nuisance parameters are suppressed for notational convenience. It is taken from Singh, Xie and Strawderman [25, 26]. The CD definition is essentially the same as in Schweder and Hjort [22]; they did not define the asymptotic CD however.

**Definition 1.1.** A function $H_n(\cdot) = H_n(\mathbf{X}_n, \cdot)$ on $\mathcal{X} \times \Theta \to [0, 1]$ is called a confidence distribution (CD) for a parameter $\theta$, if (i) For each given sample set $\mathbf{X}_n$ in the sample set space $\mathcal{X}$, $H_n(\cdot)$ is a continuous cumulative distribution function in the parameter space $\Theta$; (ii) At the true parameter value $\theta = \theta_0$, $H_n(\theta_0) = H_n(\mathbf{X}_n, \theta_0)$, as a function of the sample set $\mathbf{X}_n$, has a uniform distribution $U(0, 1)$.

The function $H_n(\cdot)$ is called an asymptotic confidence distribution (aCD), if requirement (ii) above is replaced by (ii)': At $\theta = \theta_0$, $H_n(\theta_0) \xrightarrow{W} U(0,1)$, as $n \to +\infty$, and the continuity requirement on $H_n(\cdot)$ is dropped.

We call, when it exists, $h_n(\theta) = H_n'(\theta)$ a CD density. It is also known as confidence density in the literature. It follows from the definition of CD that if $\theta < \theta_0$, $H_n(\theta) \stackrel{sto}{\leq} 1 - H_n(\theta)$, and if $\theta > \theta_0$, $1 - H_n(\theta) \stackrel{sto}{\leq} H_n(\theta)$. Here, $\stackrel{sto}{\leq}$ is a stochastic comparison between two random variables; i.e., for two random variable $Y_1$ and $Y_2$, $Y_1 \stackrel{sto}{\leq} Y_2$, if $P(Y_1 \leq t) \geq P(Y_2 \leq t)$ for all $t$. Thus a CD works, in a sense, like a compass needle. It points towards $\theta_0$, when placed at $\theta \neq \theta_0$, by assigning more mass stochastically to that side (left or right) of $\theta$ that contains $\theta_0$. When placed at $\theta_0$ itself, $H_n(\theta) = H_n(\theta_0)$ has the uniform U[0,1] distribution and thus it is noninformative in direction.

The interpretation of a CD as a distribution estimator is as follows. The purpose of analyzing sample data is to gather knowledge about the population from which the sample came. The unknown $\theta$ is a characteristic of the population. Though useful, the knowledge acquired from the data analysis is imperfect in the sense that there is still a, usually known, degree of uncertainty remaining. Statisticians can present the acquired knowledge on $\theta$, with the left-over uncertainty, in the form of a probability distribution. This appropriately calibrated distribution, that reflects statisticians' confidence regarding where $\theta$ lives, is a CD. Thus, a CD is an expression of inference (an inferential output) and not a distribution on $\theta$. What is really fascinating is that a CD is loaded with a wealth of information about $\theta$ (as it is detailed later), as is a posterior distribution in Bayesian inference.

Before we give some illustrative examples, let us describe a general substitution scheme for the construction of CDs, that avoids inversion of functions; See, also Schweder and Hjort [22]. Although this scheme does not cover all possible ways of constructing CDs (see, for example, Section 4), it covers a wide range of examples involving pivotal statistics.

Consider a statistical function $\psi(\mathbf{X}_n, \theta)$ which involves the data set $\mathbf{X}_n$ and the parameter of interest $\theta$. Besides $\theta$, the function $\psi$ may contain some known parameters (which should be treated as constants) but it should not have any other unknown parameter. On $\psi$, we impose the following condition:

- *For any given $\mathbf{X}_n$, $\psi(\mathbf{X}_n, \theta)$ is continuous and monotonic as a function of $\theta$.*

Suppose further that $G_n$, the true c.d.f. of $\psi(\mathbf{X}_n, \theta)$, does not involve any unknown parameter and it is analytically tractable. In such a case, $\psi(\mathbf{X}_n, \theta)$ is generally known as a pivot. Then one has the following exact CD for $\theta$ (provided $G_n(\cdot)$ is continuous):

$$H_n(x) = \begin{cases} G_n(\psi(\mathbf{X}_n, x)), & \text{if } \psi \text{ is increasing in } \theta \\ 1 - G_n(\psi(\mathbf{X}_n, x)), & \text{if } \psi \text{ is decreasing in } \theta. \end{cases}$$

In most cases, $H_n(x)$ is typically a continuous c.d.f. for fixed $\mathbf{X}_n$ and, as a function of $\mathbf{X}_n$, $H_n(\theta_0)$ follows a $U[0,1]$ distribution. Thus, $H_n$ is a CD by definition. Note the substitution of $\theta$ by $x$.

In case the sampling distribution $G_n$ is unavailable, including the case in which $G_n$ depends on unknown nuisance parameters, one can turn to an approximate or estimated sampling distribution $\hat{G}_n$. This could be the limit of $G_n$, an estimate of the limit or an estimate based on bootstrap or some other method. Utilizing $\hat{G}_n$, one defines

$$H_n(x) = \begin{cases} \hat{G}_n(\psi(\mathbf{X}_n, x)), & \text{if } \psi \text{ is increasing in } \theta, \\ 1 - \hat{G}_n(\psi(\mathbf{X}_n, x)), & \text{if } \psi \text{ is decreasing in } \theta. \end{cases}$$

In most cases, $H_n(\theta_0) \xrightarrow{L} U[0,1]$ and is thus an asymptotic CD. The above construction resembles Beran's construction of prepivot (see Beran [4], page 459), which was defined at $\theta = \theta_0$ (the true value of $\theta$). Beran's goal was to achieve second order accuracy in general via double bootstrap.

We now present some illustrative examples of CDs.

**Example 1.1.** (*Normal mean and variance*) The most basic case is that of sampling from a normal distribution with parameters $\mu$ and $\sigma^2$. Consider first a CD for the mean when the variance is unknown. Here the standard pivot is $\psi(\mathbf{X}_n, \mu) = (\bar{X}_n - \mu)/(s_n/\sqrt{n})$, which has the student $t$-distribution with $(n-1)$ d.f. Using the above substitution, the CD for $\mu$ is

$$H_n(x) = 1 - P\left(T_{n-1} \leq \frac{\bar{X} - x}{s_n/\sqrt{n}}\right) = P\left(T_{n-1} \leq \frac{x - \bar{X}}{s_n/\sqrt{n}}\right),$$

where $T_{n-1}$ is a random variable that has the Student's $t_{n-1}$- distribution.

For $\sigma^2$, the usual pivot is $\psi(\mathbf{X}_n, \sigma^2) = (n-1)s_n^2/\sigma^2$. By the substitution method, the CD for $\sigma^2$ is

(1.1) $$H_n(x) = P\left(\chi_{n-1}^2 \geq \frac{(n-1)s_n^2}{x}\right), \quad x \geq 0.$$

where $\chi_{n-1}^2$ is a random variable that has the Chi-square distribution with $n-1$ degrees of freedom.

**Example 1.2.** (*Bivariate normal correlation*) For a bivariate normal population, let $\rho$ denote the correlation coefficient. The asymptotic pivot used in this example is Fisher's Z, $\psi(\mathbf{X}_n, \rho) = \left[\log((1+r)/(1-r)) - \log{(1+\rho)/(1-\rho)}\right]/2$ where $r$ is the sample correlation. Its limiting distribution is $N\left(0, \frac{1}{n-3}\right)$, with a fast rate of convergence. So the resulting asymptotic CD is

$$H_n(x) = 1 - \Phi\left(\sqrt{n-3}\left(\frac{1}{2}\log\frac{1+r}{1-r} - \frac{1}{2}\log\frac{1+x}{1-x}\right)\right), \quad -1 \leq x \leq 1.$$

**Example 1.3.** (*Nonparametric bootstrap*) Turning to nonparametric examples based on bootstrap, let $\hat{\theta}$ be an estimator of $\theta$, such that the limiting distribution of $\hat{\theta}$, properly normalized, is symmetric. Using symmetry, if the sampling distribution of $\hat{\theta} - \theta$ is estimated by the bootstrap distribution of $\hat{\theta} - \hat{\theta}_B$, then an asymptotic CD is given by

$$H_n(x) = 1 - P_B(\hat{\theta} - \hat{\theta}_B \leq \hat{\theta} - x) = P_B(\hat{\theta}_B \leq x).$$

Here, $\hat{\theta}_B$ is $\hat{\theta}$ computed on a bootstrap sample. The resulting asymptotic CD is the raw bootstrap distribution of $\hat{\theta}$.

If the distribution of $\hat{\theta} - \theta$ is estimated by the bootstrap distribution of $\hat{\theta}_B - \hat{\theta}$, which is what bootstrappers usually do, the corresponding asymptotic CD is

$$H_n(x) = 1 - P_B(\hat{\theta}_B - \hat{\theta} \leq \hat{\theta} - x) = P_B(\hat{\theta}_B \geq 2\hat{\theta} - x).$$

**Example 1.4.** (*Bootstrap-t method*) By the bootstrap-$t$ method, the distribution of asymptotic pivot $(\hat{\theta} - \theta)/\widehat{SE}(\hat{\theta})$ is estimated by the bootstrap distribution of $(\hat{\theta}_B - \hat{\theta})/\widehat{SE}_B(\hat{\theta}_B)$. Here $\widehat{SE}_B(\hat{\theta}_B)$ is the estimated standard error of $\hat{\theta}_B$, based on the bootstrap sample. Such an approximation has so-called, second order accuracy (see Singh [24], Babu and Singh [2, 3]). The resulting asymptotic CD would be

$$H_n(x) = 1 - P_B\left(\frac{\hat{\theta}_B - \hat{\theta}}{\widehat{SE}_B(\hat{\theta}_B)} \leq \frac{\hat{\theta} - x}{\widehat{SE}(\hat{\theta})}\right).$$

Such a CD, at $x = \theta_0$, typically converges to $U[0, 1]$, in law, at a rapid pace.

**Example 1.5.** (*Bootstrap 3rd order accurate aCD*) Hall [13] came up with the following increasing function of the $t$-statistics, which does not have the $1/\sqrt{n}$-term in its Edgeworth expansion:

$$\psi(\mathbf{X}_n, \mu) = t + \frac{\hat{\lambda}}{6\sqrt{n}}(2t^2 + 1) + \frac{1}{27n}\hat{\lambda}^2 t^3.$$

Here $t = \sqrt{n}(\bar{X} - \mu)/s_n$, $\lambda = \frac{\mu_3}{\sigma^3}$, $\hat{\lambda}$ is a sample estimate of $\lambda$ and the assumption of population normality is dropped. Under mild conditions on the population distribution, the bootstrap approximation to the distribution of this function of $t$, is third-order correct. Let $\hat{G}_B$ be the c.d.f. of the bootstrap approximation. Then, using the substitution, a second-order correct CD for $\mu$ is given by

$$H_n(x) = 1 - \hat{G}_B(\psi(\mathbf{X}_n, x)).$$

One also has CDs that do not involve pivotal statistics. A particular class of such CDs are constructed from likelihood functions. We will have some detailed discussions on the connections of CDs and likelihood functions in Section 4.

For each given sample $\mathbf{X}_n$, $H_n(\cdot)$ is a cumulative distribution function. We can construct a random variable $\xi$ such that $\xi$ has the distribution $H_n$. For convenience of presentations, we call $\xi$ a CD random variable.

**Definition 1.2.** We call $\xi = \xi_{H_n}$ a CD random variable associated with a CD $H_n$, if the conditional distribution of $\xi$ given the data $\mathbf{X}_n$ is $H_n$.

As an example, let $U$ be a $U(0,1)$ random variable independent of $\mathbf{X}_n$, then $\xi = H_n^{-1}(U)$ is a CD random variable.

Let us note that $\xi$ may be viewed as a CD-randomized estimator of $\theta_0$. As an estimator, $\xi$ is median unbiased, i.e., $P_{\theta_0}(\xi \leq \theta_0) = E_{\theta_0}\{H_n(\theta_0)\} = \frac{1}{2}$. However, $\xi$ is not always mean unbiased. For example, the CD random variable $\xi$ associated with (1.1) in Example 1.1 is mean biased as an estimator of $\sigma^2$.

We close this section with a equivariance result on CDs, which may be helpful in the construction of a CD for a function of $\theta$. For example, to derive a CD for $\sigma$ from that of $\sigma^2$ given in Example 1.1. The equivariance is shared by Efron's bootstrap distribution of an estimator, which is of course aCD (Example 1.3 above) under conditions.

**Proposition 1.1.** *Let $H_n$ be a CD for $\theta$ and $\xi$ be an associated CD random variable. Then, the conditional distribution function of $g(\xi)$, for given $\mathbf{X}_n$, is a CD of $g(\theta)$, if $g$ is monotonic. When the monotonicity is limited to a neighborhood of $\theta_0$ only, then the conditional distribution of $g(\xi)$, for given $\mathbf{X}_n$, yields an asymptotic CD at $\theta = \theta_0$, provided, for all $\epsilon > 0$, $H_n(\theta_0 + \epsilon) - H_n(\theta_0 - \epsilon) \xrightarrow{p} 1$.*

*Proof.* The proof of the first claim is straightforward. For the second claim, we note that, if $g(\cdot)$ is increasing within $(\theta_0 - \epsilon, \theta_0 + \epsilon)$, $P(g(\xi) \leq g(\theta_0)|\mathbf{x}) = P(\{\xi \leq \theta_0\} \cap \{\theta_0 - \epsilon \leq \xi \leq \theta_0 + \epsilon\}|\mathbf{x}) + o_p(1) = H_n(\theta_0) + o_p(1)$. One argues similarly for decreasing $g(\cdot)$. □

The rest of the paper is arranged as follows. Section 2 is devoted to comparing CDs for the same parameter and related issues. In Section 3, we explore, from the frequentist viewpoint, inferential information contained within a CD. In Section 4, we establish that the normalized profile likelihood function is an aCD. Lastly, Section 5 is an attempt to formally define and develop the notion of joint CD for a parameter vector. Parts of Sections 2 and 3 are closely related to the recent paper of Schweder and Hjort [22], and also to Singh, Xie and Strawderman [25, 26]. Schweder and Hjort [22] present essentially the same definition of the CD and also compare CDs as we do in this paper (See Definition 2.1). They also develop the notion of an optimal CD which is quite close to that presented here and in Singh, Xie and Strawderman [25]. Our development is based on the theory of UMPU tests and differs slightly from theirs. The materials on p-values in Section 3.3 is also closely related to, but somewhat more general than, that of Fraser [11].

## 2. Comparison of CDs and a notion of optimal CD

The precision of a CD can be measured in terms of how little probability mass a CD wastes on sets that do not include $\theta_0$. This suggests that, for $\epsilon > 0$, one should compare the quantities $H_1(\theta_0 - \epsilon)$ with $H_2(\theta_0 - \epsilon)$ and also $1 - H_1(\theta_0 + \epsilon)$ and $1 - H_2(\theta_0 + \epsilon)$. In each case, a smaller value is preferred. Here $H_1$ and $H_2$ are any two CDs for the common parameter $\theta$, based on the same sample of size $n$.

**Definition 2.1.** Given two CDs $H_1$ and $H_2$ for $\theta$, we say $H_1$ is more precise than $H_2$, at $\theta = \theta_0$, if for all $\epsilon > 0$,

$$H_1(\theta_0 - \epsilon) \stackrel{\text{sto}}{\leq} H_2(\theta_0 - \epsilon) \quad \text{and} \quad 1 - H_1(\theta_0 + \epsilon) \stackrel{\text{sto}}{\leq} 1 - H_2(\theta_0 + \epsilon)$$

when $\theta_0$ is the prevailing value of $\theta$.

An essentially equivalent definition is also used in Singh, Xie and Strawderman [25] and Schweder and Hjort [22]. The following proposition follows immediately from the definition.

**Proposition 2.1.** *If $H_1$ is more precise than $H_2$ and they both are strictly increasing, then for all $t$ in $[0, 1]$,*

$$\left[H_1^{-1}(t) - \theta_0\right]^+ \overset{sto}{\leq} \left[H_2^{-1}(t) - \theta_0\right]^+ \text{ and } \left[H_1^{-1}(t) - \theta_0\right]^- \overset{sto}{\geq} \left[H_2^{-1}(t) - \theta_0\right]^-$$

*Thus,*

$$(2.1) \qquad \left|H_1^{-1}(t) - \theta_0\right| \overset{sto}{\leq} \left|H_2^{-1}(t) - \theta_0\right|.$$

The statement (2.1) yields a comparison of confidence intervals based on $H_1$ and $H_2$. In general, an endpoint of a confidence interval based on $H_1$ is closer to $\theta_0$ than that based on $H_2$. In some sense, it implies that the $H_1$ based confidence intervals are more compact. Also, the CD median of a more precise CD is stochastically closer to $\theta_0$.

Let $\phi(x, \theta)$ be a loss function such that $\phi(\cdot, \cdot)$ is non-decreasing for $x \geq \theta$ and non-increasing for $x \leq \theta$. We now connect the above defined CD-comparison to the following concept of the $\phi$-dispersion of a CD.

**Definition 2.2.** *For a CD $H(x)$ of a parameter $\theta$, the $\phi$-dispersion of $H(x)$ is defined as*

$$d_\phi(\theta, H) = E_\theta \int \phi(x, \theta) dH(x).$$

In the special case of square error loss, $d_{sq}(\theta, H) = E_\theta \int (x-\theta)^2 dH(x)$. In general, we have the following:

**Theorem 2.1.** *If $H_1$ is more precise than $H_2$ at $\theta = \theta_0$, in terms of Definition 2.1, then*

$$(2.2) \qquad d_\phi(\theta_0, H_1) \leq d_\phi(\theta_0, H_2).$$

In fact, the above theorem holds under a set of weaker conditions: For any $\epsilon > 0$,

$$(2.3) \quad E\{H_1(\theta_0 - \epsilon)\} \leq E\{H_2(\theta_0 - \epsilon)\} \text{ and } E\{H_1(\theta_0 + \epsilon)\} \geq E\{H_2(\theta_0 + \epsilon)\}$$

*Proof.* The claim in (2.2) is equivalent to

$$(2.4) \qquad E\{\phi(\xi_1, \theta_0)\} \leq E\{\phi(\xi_2, \theta_0)\},$$

where $\xi_1$ and $\xi_2$ are CD random variables associated with $H_1$ and $H_2$ (see, Definition 1.2), respectively. From (2.3), via conditioning on $\mathbf{X}_n$, it follows that

$$(\xi_1 - \theta_0)^+ \overset{sto}{\leq} (\xi_2 - \theta_0)^+ \text{ and } (\xi_1 - \theta_0)^- \overset{sto}{\geq} (\xi_2 - \theta_0)^-.$$

Due to the monotonicity of $\phi(\cdot, \theta_0)$, we have

$$\phi(\xi_1, \theta_0) I_{(\xi_1 \geq \theta_0)} \overset{sto}{\leq} \phi(\xi_2, \theta_0) I_{(\xi_2 \geq \theta_0)} \text{ and } \phi(\xi_1, \theta_0) I_{(\xi_1 < \theta_0)} \overset{sto}{\leq} \phi(\xi_2, \theta_0) I_{(\xi_2 < \theta_0)}.$$

The above inequalities lead to (2.4) immediately. □

Suppose now that there is a family of Uniformly Most Powerful Unbiased tests for testing $K_0 : \theta \leq \theta_0$ versus $K_1 : \theta > \theta_0$, for every $\theta_0$. The underlying family of distributions may have nuisance parameter(s). Let the corresponding $p$-value (the inf of $\alpha$ at which $K_0$ can be rejected) $p(\theta_0) = p(\mathbf{X}_n, \theta_0)$ be strictly increasing and continuous as a function of $\theta_0$. It is further assumed that $1 - p(\theta_0)$ is the $p$-value of an UMPU test for testing $K_0 : \theta \geq \theta_0$ vs $K_1 : \theta < \theta_0$. Let the distribution of $p(\theta_0)$ under $\theta_0$ be $U[0,1]$ and let the range of $p(\cdot)$ be $[0,1]$. Define the corresponding CD, $H^*(x) = p(\mathbf{X}_n, x)$. We have the following result.

**Theorem 2.2.** *The CD $H^*$ defined above is more precise than any other CD for the parameter $\theta$, at all $\theta_0$.*

*Proof.* Let $\theta = \theta_0$ be the true value. Note that $P_{\theta_0}\big(H^*(\theta_0 - \epsilon) < \alpha\big)$ is the power (at $\theta = \theta_0$) of the UMPU test when $K_0$ is $\theta \leq \theta_0 - \epsilon$ and $K_1$ is $\theta > \theta_0 - \epsilon$. Given any other CD $H$, one has the following unbiased test for testing the same hypotheses: Reject $K_0$ iff $H(\theta_0 - \epsilon) < \alpha$. Therefore, $P_{\theta_0}\big(H^*(\theta_0 - \epsilon) < \alpha\big) \geq P_{\theta_0}\big(H(\theta_0 - \epsilon) < \alpha\big)$ for all $\alpha \in [0,1]$. Using the function $1 - p(\cdot)$, one similarly argues for $P_{\theta_0}\big(1 - H^*(\theta_0 + \epsilon) < \alpha\big) \geq P_{\theta_0}\big(1 - H(\theta_0 + \epsilon) < \alpha\big)$. Thus, $H^*$ is most precise. □

It should be mentioned that the property of CDs as exhibited in Theorem 2.2 depend on corresponding optimality properties of hypothesis tests. The basic ideas behind this segment could be traced to the discussions of confidence intervals in Lehmann [14].

**Remark 2.1.** If the underlying parametric family has the so-called MLR (monotone likelihood ratio) property, there exists an UMP test for one-sided hypotheses whose $p$-value is monotonic.

**Example 2.1.** In the testing problem of normal means, the $Z$-test is UMPU (actually UMP), for the one-sided hypotheses when $\sigma$ is known. The $t$-test is UMPU for the one-sided hypotheses when $\sigma$ is a nuisance parameter (see Lehmann [14], Chapter 5). The conclusion: $H^*(x) = \Phi\big(\frac{x - \bar{X}}{\sigma/\sqrt{n}}\big)$ is the most precise CD for $\mu$, when $\sigma$ is known, and $H^{**}(x) = F_{t_{n-1}}\big(\frac{x - \bar{X}}{s_n/\sqrt{n}}\big)$ is the most precise CD, when $\sigma$ is not known. Here, $F_{t_{n-1}}$ is the cumulative distribution function of the $t$-distribution with degrees of freedom $n - 1$.

The above presented optimality theory can be expressed in the decision theoretic framework by considering the "target distribution" towards which a CD is supposed to converge. Given $\theta_0$ as the true value of $\theta$, the target distribution is $\delta(\theta_0)$, the Dirac $\delta$-measure at $\theta_0$, which assigns its 100% probability mass at $\theta_0$ itself. A loss function can be defined in terms of "distance" between a CD $H(\cdot)$ and $\delta(\theta_0)$.

Perhaps, the most popular distance, between two distributions $F$ and $G$, is the Kolmogorov-Smirnov distance $= \sup_x |F(x) - G(x)|$. However, it turns out that this particular distance, between a CD and its target $\delta(\theta_0)$, is useless for comparing two CDs. To see this, note that

$$E_{\theta_0}\left[\sup_x |H(x) - I_{[\theta_0, \infty)}|\right] = E_{\theta_0}\big[\max\big(H(\theta_0), 1 - H(\theta_0)\big)\big] = 3/4 \quad \text{(free of } H!\text{)},$$

since $H(\theta_0)$ follows the $U[0,1]$ distribution. Note, $I_{[\theta_0, \infty)}$ is the cdf of $\delta(\theta_0)$. So, we instead consider the integrated distance

$$\tau\big(H, \delta(\theta_0)\big) = \int \psi\big(|H(x) - I_{[\theta_0, \infty)}|\big)\, dW(x)$$

where $\psi(\cdot)$ is a monotonic function from $[0,1]$ to $R^+$ and $W(\cdot)$ is a positive measure. The risk function is

$$R_{\theta_0}(H) = E_{\theta_0}[\tau(H, \delta(\theta_0))].$$

**Theorem 2.3.** *If $H_1$ is a more precise CD than $H_2$ at $\theta = \theta_0$, then $R_{\theta_0}(H_1) \leq R_{\theta_0}(H_2)$.*

Theorem 2.3 is proved by interchanging the expectation and the integration appearing in the loss function, which is allowed by Fubini's Theorem.

Now, for asymptotic CDs, we define two asymptotic notions of "more precise CD". One is in the Pitman sense of "local alternatives" and the other is in the Bahadur sense when the $\theta \neq \theta_0$ is held fixed and the a.s. limit is taken on CDs themselves.

First, the Pitman-more precise CD:

**Definition 2.3.** Let $H_1$ and $H_2$ be two asymptotic CDs, We say that $H_1$ is Pitman-more precise than $H_2$ if, for every $\epsilon > 0$,

$$\lim P\left(H_{1n}\left(\theta_0 - \frac{\epsilon}{\sqrt{n}}\right) \leq t\right) \geq \lim P\left(H_{2n}\left(\theta_0 - \frac{\epsilon}{\sqrt{n}}\right) \leq t\right)$$

and

$$\lim P\left(1 - H_{1n}\left(\theta_0 + \frac{\epsilon}{\sqrt{n}}\right) \leq t\right) \geq \lim P\left(1 - H_{2n}\left(\theta_0 + \frac{\epsilon}{\sqrt{n}}\right) \leq t\right)$$

where all the limits (as $n \to \infty$) are assumed to exist, and the probabilities are under $\theta = \theta_0$.

Thus, we are requiring that in terms of the limiting distributions,

$$H_{1n}\left(\theta_0 - \frac{\epsilon}{\sqrt{n}}\right) \overset{sto}{\leq} H_{2n}\left(\theta_0 - \frac{\epsilon}{\sqrt{n}}\right)$$

and

$$1 - H_{1n}\left(\theta_0 + \frac{\epsilon}{\sqrt{n}}\right) \overset{sto}{\leq} 1 - H_{2n}\left(\theta_0 + \frac{\epsilon}{\sqrt{n}}\right).$$

**Example 2.2.** The most basic example allowing such a comparison is that of

$$H_{1n}(x) = \Phi\left(\frac{x - \hat{\theta}_1}{a/\sqrt{n}}\right) \quad \text{and} \quad H_{2n}(x) = \Phi\left(\frac{x - \hat{\theta}_2}{b/\sqrt{n}}\right)$$

where $\hat{\theta}_1, \hat{\theta}_2$ are two $\sqrt{n}$-consistent asymptotically normal estimators of $\theta$, with asymptotic variances $a^2/n$ and $b^2/n$, respectively. The one with a smaller asymptotic variance is Pitman-more precise.

Next, the Bahadur-type comparison:

**Definition 2.4.** We define $H_1$ to be Bahadur-more precise than $H_2$ (when $\theta = \theta_0$) if, for every $\epsilon > 0$, a.s.

$$\lim \frac{1}{n} \log H_1(\theta_0 - \epsilon) \leq \lim \frac{1}{n} \log H_2(\theta_0 - \epsilon)$$

and

$$\lim \frac{1}{n} \log \left(1 - H_1(\theta_0 + \epsilon)\right) \leq \lim \frac{1}{n} \log \left(1 - H_2(\theta_0 + \epsilon)\right),$$

where the limits, as $n \to \infty$, are assumed to exist.

Here too, we are saying that in a limit sense, $H_1$ places less mass than $H_2$ does on half-lines which exclude $\theta_0$. This comparison is on $H_i$ directly and not on their distributions.

**Example 2.3.** Let us return to the example of $Z$ vs $t$ for normal means. The fact that $\Phi\left(\frac{x-\bar{X}}{\sigma/\sqrt{n}}\right)$ is Bahadur-more precise than $F_{t_{n-1}}\left(\frac{x-\bar{X}}{s_n/\sqrt{n}}\right)$ follows from the well-known limits:

$$\frac{1}{n}\log \Phi(-\epsilon\sqrt{n}) \longrightarrow -\epsilon^2/2 \quad \text{and} \quad \frac{1}{n}\log F_{t_{n-1}}(-\epsilon\sqrt{n}) \longrightarrow -\frac{1}{2}\log\left(1+\epsilon^2\right).$$

**Remark 2.2.** Under modest regularity conditions,

$$\frac{1}{n}\lim \log H_n(\theta_0 - \epsilon) = \frac{1}{n}\lim \log h_n(\theta_0 - \epsilon)$$

and

$$\frac{1}{n}\lim \log[1 - H_n(\theta_0 + \epsilon)] = \frac{1}{n}\lim \log h_n(\theta_0 + \epsilon).$$

The right hand sides are $CD$-density slopes, which have significance of their own for CDs. The faster the CD density goes to 0, at fixed $\theta \neq \theta_0$, the more compact, in limit, the CD is.

## 3. Information contained in a CD

This section discusses inference on $\theta$ from a CD or an aCD. We briefly consider basic elements of inferences about $\theta$, including confidence intervals, point estimation, and hypothesis testing.

### 3.1. Confidence intervals

The derivation of confidence intervals from a CD is straightforward and well known. Note that, according to the definition, the intervals $(-\infty, H_n^{-1}(\alpha)]$ and $[H_n^{-1}(\alpha), +\infty)$ are one-sided confidence intervals for $\theta$, for any $\alpha \in (0,1)$. It is also clear that the central regions of the CD, i.e., $(H_n^{-1}(\alpha/2), H_n^{-1}(1 - \alpha/2))$, provide two sided confidence interval for $\theta$ at each coverage level $\alpha \in (0,1)$. The same is true for an aCD, where the confidence level is achieved in limit.

### 3.2. Point estimators

A CD (or an aCD) on the real line can be a tool for point estimation as well. We assume the following condition, which is mild and almost always met in practice.
(3.1) *For any $\epsilon$ and each fixed $\theta_0$, $0 < \epsilon < \frac{1}{2}$, $L_n(\epsilon) = H_n^{-1}(1 - \epsilon) - H_n^{-1}(\epsilon) \to 0$, in probability, as $n \to \infty$.*
Condition (3.1) states the CD based information concentrates around $\theta_0$ as $n$ gets large.

One natural choice for a point estimator of the parameter $\theta$ is the median of a CD (or an aCD), $M_n = H_n^{-1}(1/2)$. Note that, $M_n$ is a median-unbiased estimator; even if the original estimator, on which $H_n$ is based, is not. For instance, this does happen in the case of the CD for the normal variance, based on the $\chi^2$-distribution. This median unbiased result follows from observation that $P_{\theta_0}(M_n \leq \theta_0) = P_{\theta_0}(1/2 \leq H_n(\theta_0)) = 1/2$. The following is a consistency result on $M_n$.

**Theorem 3.1.** *(i) If condition (3.1) is true, then $M_n \to \theta_0$, in probability, as $n \to \infty$. (ii) Furthermore, if $L_n(\epsilon) = O_p(a_n)$, for a non-negative $a_n \to 0$, then $M_n - \theta_0 = O_p(a_n)$.*

*Proof.* (i) We first note the identity: for $\alpha \in (0, \frac{1}{2})$

$$P_{\theta_0}(|M_n - \theta_0| > \delta) = P_{\theta_0}\big(\{|M_n - \theta_0| > \delta\} \cap \{L_n(\alpha) > \delta\}\big) \\ + P_{\theta_0}\big(\{|M_n - \theta_0| > \delta\} \cap \{L_n(\alpha) \leq \delta\}\big).$$

Under the assumed condition, the first term in the r.h.s.$\to 0$. We prove that the second term is $\leq 2\alpha$. This is deduced from the set inequality:

$$\{|M_n - \theta_0| > \delta\} \cap \{L_n(\alpha) \leq \delta\} \subseteq \{H_n(\theta_0) \leq \alpha\} \bigcup \{H_n(\theta_0) \geq 1 - \alpha\}.$$

To conclude the above set inequality, one needs to consider the two cases: $M_n > \theta_0 + \delta$ and $M_n < \theta_0 - \delta$, separately. The first case leads to $\{H_n(\theta_0) \leq \alpha\}$ and the second one to $\{H_n(\theta_0) \geq 1 - \alpha\}$. Part (i) follows, since $\alpha$ is arbitrary.

One can prove part (ii) by using similar reasoning. □

One can also use the average of a CD (or an aCD), $\bar{\theta}_n = \int_{-\infty}^{+\infty} t \, dH_n(t)$, to consistently estimate the unknown $\theta_0$. Indeed, $\bar{\theta}_n$ is the frequentist analog of Bayesian estimator of $\theta$ under the usual square loss.

**Theorem 3.2.** *Under condition (3.1), if $r_n = \int_{-\infty}^{+\infty} t^2 dH_n(t)$ is bounded in probability, then $\bar{\theta}_n \to \theta_0$, in probability.*

*Proof.* Using Cauchy Schwartz inequality, we have, for any $0 < \epsilon < 1/2$,

$$\big| \int_{-\infty}^{H_n^{-1}(\epsilon)} t dH_n(t) + \int_{H_n^{-1}(1-\epsilon)}^{+\infty} t dH_n(t) \big| \leq 2 r_n \epsilon^{1/2}.$$

Thus,

$$(1 - \epsilon) H_n^{-1}(\epsilon) - 2 r_n \epsilon^{1/2} \leq \bar{\theta}_n \leq (1 - \epsilon) H_n^{-1}(1 - \epsilon) + 2 r_n \epsilon^{1/2}.$$

Now, with $M_n = H_n^{-1}(1/2)$,

$$|\bar{\theta}_n - \theta_0| \leq |M_n - \theta_0| + |\bar{\theta}_n - M_n| \leq |M_n - \theta_0| + |H_n^{-1}(1-\epsilon) - H_n^{-1}(\epsilon)| + 2 r_n \epsilon^{1/2} + 2\epsilon |M_n|$$

Since $\epsilon > 0$ is arbitrary, the result follows using Theorem 3.1. □

Denote $\hat{\theta}_n = \arg\max_\theta h_n(\theta)$, the value that maximizes the CD (or aCD) density function $h_n(\theta) = \frac{d}{d\theta} H_n(\theta)$. Let $\epsilon_n = \inf_{0 < \epsilon \leq 1/2} \{\epsilon : \hat{\theta}_n \notin [H_n^{-1}(\epsilon), H_n^{-1}(1-\epsilon)]\}$. The event $\epsilon_n > \epsilon^*$ is that $\hat{\theta}_n$ will not be in the tails having probability less than $\epsilon^*$. We have the following theorem.

**Theorem 3.3.** *Assume condition (3.1) holds. Suppose there exists a fixed $\epsilon^* > 0$, such that $P(\epsilon_n \geq \epsilon^*) \to 1$. Then, $\hat{\theta}_n \to \theta_0$ in probability.*

*Proof.* Note that $\hat{\theta}_n \in [H_n^{-1}(\epsilon^*), H_n^{-1}(1 - \epsilon^*)]^c$ implies $\epsilon_n \leq \epsilon^*$. The claim follows immediately using (3.1) and Theorem 3.1. □

### 3.3. Hypothesis testing

Now, let us turn to one-sample hypothesis testing, given that a CD $H_n(\cdot)$ is available on a parameter of interest $\theta$. Suppose the null hypothesis is $K_0 : \theta \in C$ versus the alternative $K_1 : \theta \in C^c$. A natural line of thinking would be to measure the support that $H_n(\cdot)$ lends to $C$. If the support is "high," the verdict based on $H_n$ should be for $C$ and if it is low, it should be for $C^c$. The following two definitions of support for $K_0$ from a CD are suggested by classical p-values. These two notions of support highlight the distinction between the two kind of p-values used in statistical practice, one for the one-sided hypotheses and the other for the point hypotheses.

I. *Strong-support $p_s(C)$*, defined as $p_s(C) = H_n(C)$, which the probability content of $C$ under $H_n$.

II. *Weak-support $p_w(C)$*, defined as $p_w(C) = \sup_{\theta \in C} 2\min(H_n(\theta), 1 - H_n(\theta))$.

See, e.g., Cox and Hinkley [6], Barndorff-Nielsen and Cox [1], and especially Fraser [11], for discussions on this topic related to p-value functions. Our results are closely related to those in Fraser [11] but they are developed under a more general setting. We use the following claim for making connection between the concepts of support and the p-values.

**Claim A** If $K_0$ is of the type $(-\infty, \theta_0]$ or $[\theta_0, \infty)$, the classical p-value typically agrees with the strong-support $p_s(C)$. If $K_0$ is a singleton, i.e. $K_0$ is $\theta = \theta_0$, then the p-value typically agrees with the weak-support $p_w(C)$.

To illustrate the above claim, consider tests based on the normalized statistic $T_n = (\hat{\theta} - \theta)/SE(\hat{\theta})$, for an arbitrary estimator $\hat{\theta}$. Based on the method given in Section 1.1, a CD for $\theta$ is $H_n(x) = P_{\eta_n}(\hat{\theta} - SE(\hat{\theta})\eta_n \leq x)$, where $\eta_n$ is independent of the data $\mathbf{X}_n$ and $\eta_n \overset{\mathcal{L}}{=} T_n$. Thus,

$$p_s(-\infty, \theta_0) = H_n(\theta_0) = P_{\eta_n}\left(\eta_n \geq \frac{\hat{\theta} - \theta_0}{SE(\hat{\theta})}\right).$$

This agrees with the p-values for one-sided test $K_0 : \theta \leq \theta_0$ versus $K_1 : \theta > \theta_0$. Similar demonstrations can be given for the tests based on studentized statistics. If the null hypothesis is $K_0 : \theta = \theta_0$ vs $K_1 : \theta \neq \theta_0$, the standard p-values, based on $T_n$, is twice the tail probability beyond $(\hat{\theta} - \theta_0)/SE(\hat{\theta})$ under the distribution of $T_n$. This equals

$$2\min\left[H_n(\theta_0), 1 - H_n(\theta_0)\right] = p_w(\theta_0).$$

**Remark 3.1.** It should be remarked here that $p_w(\theta_0)$ is the Tukey's depth (see Tukey [27]) of the point $\theta_0$ w.r.t. $H_n(\cdot)$.

The following inequality justifies the names of the two supports,

**Theorem 3.4.** *For any set $C$, $p_s(C) \leq p_w(C)$.*

*Proof.* Suppose the sup in the definition of $p_w(C)$ is attained at $\theta'$, which may or may not be in $C$ and $\theta' \leq M_n$ (recall $M_n$ = the median of $H_n$). Let $\theta''$ be the point to the right of $M_n$ such that $1 - H_n(\theta'') = H_n(\theta')$. Then $C \subseteq (-\infty, \theta'] \bigcup [\theta'', \infty)$; ignoring a possible null set under $H_n$. As a consequence

$$p_s(C) = H_n(C) \leq H_n((-\infty, \theta']) + H_n([\theta'', \infty)) = p_w(C).$$

Similar arguments are given when $\theta' \geq M_n$. □

The next three theorems justify Claim A.

**Theorem 3.5.** *Let $C$ be of the type $(-\infty, \theta_0]$ or $[\theta_0, \infty)$. Then $\sup_{\theta \in C} P_\theta(p_s(C) \leq \alpha) = \alpha$.*

*Proof.* Let $C = (-\infty, \theta_0]$. For a $\theta \leq \theta_0$, $P_\theta\big(H_n((-\infty, \theta_0]) \leq \alpha\big) \leq P_\theta\big(H_n((-\infty, \theta]) \leq \alpha\big) = \alpha$. When $\theta = \theta_0$, one has the equality. A similar proof is given when $C = [\theta_0, \infty)$. □

A limiting result of the same nature holds for a more general null hypothesis, namely a union of finitely many disjoint closed intervals (bounded or unbounded). Assume the following regularity condition: as $n \to \infty$,

$$(3.2) \qquad \sup_{\theta \in [a,b]} P_\theta\big(\max\{H_n(\theta - \epsilon), 1 - H_n(\theta + \epsilon)\} > \delta\big) \to 0$$

for any finite $a, b$ and positive $\epsilon, \delta$.

Essentially, the condition (3.2) assumes that the scale of $H_n(\cdot)$ shrinks to 0, uniformly in $\theta$ lying in a compact set.

**Theorem 3.6.** *Let $C = \bigcup_{j=1}^k I_j$ where $I_j$ are disjoint intervals of the type $(-\infty, a]$ or $[c, d]$ or $[b, \infty)$. If the regularity condition (3.2) holds, then $\sup_{\theta \in C} P_\theta(p_s(C) \leq \alpha) \to \alpha$, as $n \to \infty$.*

*Proof.* It suffices to prove the claim with the sup over $\theta \in I_j$ for each $j = 1, \ldots, k$. Consider first the case when $I_j = (-\infty, a]$. For this $I_j$ and any $\delta > 0$,

$$\sup_{\theta \in I_j} P_\theta(p_s(C) \leq t) \geq P_{\theta=a}(p_s(C) \leq t) \geq P_{\theta=a}(p_s(I_j) \leq t - \delta) + o(1) = t - \delta + o(1).$$

The second inequality is due to (3.2). Also, from Theorem 3.5,

$$\sup_{\theta \in I_j} P_\theta(p_s(C) \leq t) \leq \sup_{\theta \in I_j} P_\theta(p_s(I_j) \leq t) = t$$

which completes the proof for this $I_j$. The case of $I_j = [b, \infty)$ is handled similarly. Turning to the case $I_j = [c, d]$, $c < d$, we write it as the union of $I_{j1} = \left[c, \frac{c+d}{2}\right]$ and $I_{j2} = \left[\frac{c+d}{2}, d\right]$, and note that, for any $\delta > 0$, it follows from (3.1) that

$$\sup_{\theta \in I_{j1}} P_\theta\big(p_s(C) \leq t\big) \geq P_{\theta=c}\big(p_s(C) \leq t\big) \geq P_{\theta=c}\big(p_s((c, \infty])$$

$$\leq t - \delta\big) + o(1) = t - \delta + o(1)$$

Furthermore, from Theorem 3.5, we have for any $\delta > 0$,

$$\sup_{\theta \in I_{j1}} P_\theta\big(p_s(C) \leq t\big) \leq \sup_{\theta \in I_{j1}} P_\theta\big(p_s(I_j) \leq t\big) \leq \sup_{\theta \in I_{j1}} P_\theta\big(p_s([c, \infty) \leq t + \delta\big) + o(1)$$

$$= t + \delta + o(1)$$

The case of sup over $\theta \in I_{j2}$ is dealt with in a similar way. In the arguments $\theta = c$ is replaced by $\theta = d$. □

**Remark 3.2.** The result of Theorem 2.6 still holds if $p_s(C)$ is replaced by $p_s^* = \max_{1 \leq j \leq k} p_s(I_j)$. The use of $p_s^*$ for $p$-value amounts to the so called Intersection Union Test (see, Berger [5]). Since $p_s^*$ as a $p$-value gives a larger rejection region than that by $p_s(\bigcup_1^k I_j)$ as a $p$-value, testing by $p_s^*$ will have better power, for the same asymptotic size. If the intervals $I_j$ are unbounded, it follows that

$$\sup_{\theta \in C} P_\theta(p_s(C) \leq t) \leq \sup_{\theta \in C} P_\theta(p_s^* \leq t) \leq t.$$

Moving on to the situation when $K_0$ is $\theta = \theta_0$ and $K_1$ is $\theta \neq \theta_0$, it is immediate that $p_w(\theta_0) = 2\min\{H_n(\theta_0), 1 - H_n(\theta_0)\}$ has the $U[0,1]$ distribution, since $H_n(\theta_0)$ does so. Thus $p_w(\theta_0)$ can be used like a $p$-value in the case of such a testing problem. In a more general situation when $K_0$ is $C = \{\theta_1, \theta_2, \ldots, \theta_k\}$ and $K_1$ is $C^c$, one has the following asymptotic result.

**Theorem 3.7.** *Let $K_0$ be $C = \{\theta_1, \ldots, \theta_k\}$ and $K_1$ be $C^c$. Assume that $H_n(\theta_i - \epsilon)$ and $1 - H_n(\theta_i + \epsilon) \xrightarrow{P} 0$ under $\theta = \theta_i$, for all $i = 1, 2, \ldots, k$. Then $\max_{\theta \in C} P_\theta(p_w(C) \leq \alpha) \to \alpha$, as $n \to \infty$.*

*Proof.* For simplicity, let $C = \{\theta_1, \theta_2\}$, $\theta_1 < \theta_2$. Under the condition, clearly $p_w(\theta_2) \leq 2\{1 - H_n(\theta_2)\} \xrightarrow{P} 0$, if $\theta = \theta_1$. Since $p_w(\theta_1)$ has the $U[0,1]$ distribution under $\theta = \theta_1$, it follows, using standard arguments, that, $\max\{p_w(\theta_1), p_w(\theta_2)\} \xrightarrow{\mathcal{L}} U[0,1]$ when $\theta = \theta_1$. The same holds under $\theta = \theta_2$. The result thus follows. $\square$

**Example 3.1.** (*Bio-equivalence*). An important example of the case where $C$, the null space, is a union of closed intervals is provided by the standard bioequivalence problem. In this testing problem $K_0$ is the region $\frac{\mu_1}{\mu_2} \in (-\infty, .8] \cup [1.25, \infty)$, where $\mu_1, \mu_2$ are the population means of bioavailability measures of two drugs being tested for equivalence.

**Example 3.2.** In the standard classification problem, the parameter space is divided into $k$-regions. The task is to decide which one contains the true value of parameter $\theta$. A natural (but probably over-simplified) suggestion is to compare CD contents of the $k$-regions and attach $\theta$ to the one which has got the maximum CD probability.

## 4. Profile likelihood functions and CDs

We examine here the connection between the concepts of profile likelihood function and asymptotic CD. Let $x_1, x_2, \ldots, x_n$ be independent sample draws from a parametric distribution with density $f_{\boldsymbol{\eta}}(x)$, $\boldsymbol{\eta}$ is a $p \times 1$ vector of unknown parameters. Suppose we are interested in a scalar parameter $\theta = s(\boldsymbol{\eta})$, where $s(\cdot)$ is a second-order differentiable mapping from $\mathbb{R}^p$ to $\Theta \subset \mathbb{R}$. To make an inference about $\theta$, one often obtains the log-profile likelihood function

$$\ell_n(\theta) = \sum_{i=1}^n \log f_{\hat{\boldsymbol{\eta}}(\theta)}(x_i), \qquad \text{where } \hat{\boldsymbol{\eta}}(\theta) = \underset{\{\boldsymbol{\eta}: s(\boldsymbol{\eta}) = \theta\}}{\arg\max} \sum_{i=1}^n \log f_{\boldsymbol{\eta}}(x_i).$$

Denote $\ell_n^*(\theta) = \ell_n(\theta) - \ell_n(\hat{\theta})$, where $\hat{\theta} = \arg\max_\theta \ell_n(\theta)$ is the maximum likelihood estimator of the unknown parameter $\theta$. We prove below that $e^{\ell_n^*(\theta)}$, after normalization (with respect to $\theta$, so that the area under its curve is one), is the density function of an aCD for $\theta$. The technique used to prove the main result is similar to that used in the proofs of Theorems 2 and 3 in Efron [7].

Let $i_n^{-1} = \frac{1}{n}\ell_n''(\hat{\theta})$ and $i_0^{-1}(\theta) = \lim_{n \to +\infty} \frac{1}{n}\ell_n''(\theta)$. Assume that the true value $\theta = \theta_0$ is in $\Theta^o$, the interior of $\Theta$, and $i_0^{-1}(\theta_0) > 0$. The key assumption is that

$$\sqrt{n}(\hat{\theta} - \theta_0)/\sqrt{i_n} \to N(0,1).$$

In addition to the regularity conditions that ensure the asymptotic normality of $\hat{\theta}$, we make the following three mild assumptions. They are satisfied in the cases of commonly used distributions.

(i) There exists a function $k(\theta)$, such that $\ell_n^*(\theta) \leq -nk(\theta)$ for all large $n$, a.s., and $\int_\Theta e^{-ck(\theta)}d\theta < +\infty$, for a constant $c > 0$.

(ii) There exists an $\epsilon > 0$, such that

$$c_\epsilon = \inf_{|\theta-\theta_0|<\epsilon} i_0^{-1}(\theta) > 0 \quad \text{and} \quad k_\epsilon = \inf_{|\theta-\theta_0|>\epsilon} k(\theta) > 0.$$

(iii) $\ell_n''(\theta)$ satisfies a Lipschitze condition of order 1 around $\theta_0$.

For $y \in \Theta$, write

$$H_n(y) = \frac{1}{\sqrt{2\pi i_n/n}} \int_{(-\infty,y]\cap\Theta} e^{-\frac{(\hat{\theta}-\theta)^2}{2i_n/n}} d\theta \quad \text{and} \quad G_n(y) = \frac{1}{c_n} \int_{(-\infty,y]\cap\Theta} e^{\ell_n^*(\theta)} d\theta,$$

where $c_n = \int_{-\infty}^{+\infty} e^{\ell_n^*(\theta)} d\theta$. We assume that $c_n < \infty$ for all $n$ and $\mathbf{X}_n$; condition (i) implies $c_n < \infty$ for $n$ large.

We prove the following theorem.

**Theorem 4.1.** $G_n(\theta) = H_n(\theta) + o_p(1)$ for each $\theta \in \Theta$.

*Proof.* We prove the case with $\Theta = (-\infty, +\infty)$; the other cases can be proved similarly. Let $\epsilon > 0$ be as in condition (ii). We first prove, for any fixed $s > 0$,

(4.1) $$\int_{-\infty}^{\theta_0-\epsilon} e^{\ell_n^*(\theta)} d\theta = O_p\left(\frac{1}{n^s}\right) \quad \text{and} \quad \int_{\theta_0+\epsilon}^{+\infty} e^{\ell_n^*(\theta)} d\theta = O_p\left(\frac{1}{n^s}\right).$$

Note, by condition (i) and (ii), when $n$ is large enough, we have $\ell_n^*(\theta) + s\log n < -nk(\theta) + nk_\epsilon/2 = -n(k(\theta) - k_\epsilon/2) \leq -nk_\epsilon/2$, for $|\theta - \theta_0| > \epsilon$. By Fatou's Lemma,

$$\limsup_{n\to+\infty} \int_{-\infty}^{\theta_0-\epsilon} e^{\ell_n^*(\theta)+s\log n} d\theta \leq \int_{-\infty}^{\theta_0-\epsilon} \limsup_{n\to+\infty} e^{\ell_n^*(\theta)+s\log n} d\theta$$
$$\leq \int_{-\infty}^{\theta_0-\epsilon} \lim_{n\to\infty} e^{-nk_\epsilon/2} d\theta = 0.$$

This proves the first equation in (4.1). The same is true for the second equation in (4.1).

In the case when $|\theta - \theta_0| < \epsilon$, by Taylor expansion

(4.2) $$\ell_n^*(\theta) = \frac{1}{2}\ell_n''(\tilde{\theta})(\theta - \hat{\theta})^2, \quad \text{for } \tilde{\theta} \text{ between } \theta \text{ and } \hat{\theta}.$$

From condition (ii) we have $\ell_n^*(\theta) \leq -\frac{n}{2}c_\epsilon(\theta - \hat{\theta})^2$, when $n$ is large. Thus, one can prove

(4.3) $$\sqrt{n} \int_{\theta_0-\epsilon}^{\theta_0-\frac{\epsilon}{\sqrt{n}}\log n} e^{\ell_n^*(\theta)} d\theta = o_p(1), \quad \text{and} \quad \sqrt{n} \int_{\theta_0+\frac{\epsilon}{\sqrt{n}}\log n}^{\theta_0+\epsilon} e^{\ell_n^*(\theta)} d\theta = o_p(1).$$

Now, consider the case when $|\theta - \theta_0| < \frac{\epsilon}{\sqrt{n}}\log n$. By (4.2) and condition (iii), one has

(4.4) $$\sqrt{n} \int_{\theta_0-\frac{\epsilon}{\sqrt{n}}\log n}^{\theta} e^{\ell_n^*(\theta)} d\theta = \sqrt{n} \int_{\theta_0-\frac{\epsilon}{\sqrt{n}}\log n}^{\theta} e^{-\frac{n(\theta-\hat{\theta})^2}{2i_n}} d\theta + o_p(1).$$

From (4.1), (4.3) and (4.4), it easily follows that

$$\frac{1}{\sqrt{2\pi i_n/n}} \int_{-\infty}^{\theta} e^{\ell_n^*(\theta)} d\theta = \frac{1}{\sqrt{2\pi i_n/n}} \int_{-\infty}^{\theta} e^{-\frac{n(\theta-\hat{\theta})^2}{2i_n}} d\theta + o_p(1) \qquad (4.5)$$
$$= H_n(\theta) + o_p(1),$$

for all $\theta \in (-\infty, \infty)$. Note that (4.5) implies that $c_n = \sqrt{2\pi i_n/n} + o_p(\frac{1}{\sqrt{n}})$. So (4.5) is, in fact, $G_n(\theta) = H_n(\theta) + o_p(1)$ for all $\theta \in (-\infty, \infty)$. This proves the theorem. □

**Remark 4.1.** At $\theta = \theta_0$, $H_n(\theta_0) = \Phi(\frac{\hat{\theta}-\theta_0}{\sqrt{i_n/n}}) + o_p(1) \xrightarrow{W} U(0,1)$. It follows from this theorem that $G_n(\theta_0) \xrightarrow{W} U(0,1)$, thus $G_n$ is an aCD.

**Remark 4.2.** It is well known that at the true value $\theta_0$, $-2\ell_n^*(\theta_0) = -2\{\ell_n(\theta_0) - \ell_n(\hat{\theta})\}$ is asymptotically equivalent to $n(\theta_0 - \hat{\theta})^2/i_n$; see, e.g., Murphy and van der Vaart [20]. In the proof of the above theorem, we need to extend this result to a shrinking neighborhood of $\theta_0$, and control the "tails" in our normalization of $e^{\ell_n^*(\theta)}$ (with respect to $\theta$, so that the area underneath the curve is 1). This normalization produces a proper distribution function, and Theorem 4.1 is an asymptotic result for this distribution function.

As a special case, the likelihood function in the family of one-parameter distributions (i.e., $\boldsymbol{\eta}$ is a scalar parameter) is proportional to an aCD density function. There is also a connection between the concepts of aCD and other types of likelihood functions, such as Efron's implied likelihood function, Schweder and Hjort's reduced likelihood function, etc. In fact, one can easily conclude from Theorems 1 and 2 of Efron [7] that in an exponential family, both the profile likelihood and the implied likelihood (Efron [7]) are aCD densities, after a normalization (with respect to $\theta$). Schweder and Hjort [22] proposed the reduced likelihood function, which itself is proportional to a CD density for a specially transformed parameter. See Welch and Peers [28] and Fisher [10] for earlier accounts of likelihood function based CDs in the case of single parameter families.

## 5. Multiparameter joint CD

Let us first note that in higher dimensions, the cdf is not as useful a notion, at least for our purposes here. The main reasons are: (a) The region $F(\mathbf{x}) \leq \alpha$ is not of much interest in $\mathbb{R}^k$. (b) The property $F(\mathbf{X}) \stackrel{L}{=} U[0,1]$, when $\mathbf{X} \stackrel{L}{=} F$ is lost!

### 5.1. Development of multiparameter CD through Cramer-Wold device

The following definition of a multiparameter CD has the make of a random vector having a particular multivariate distribution (arguing via characteristic functions).

**Definition 5.1.** A distribution $H_n(\cdot) \equiv H_n(\mathbf{X}_n, \cdot)$ on $\mathbb{R}^k$ is a CD in the linear sense (*l*-CD) for a $k \times 1$ parameter vector $\boldsymbol{\theta}$ if and only if for any $k \times 1$ vector $\boldsymbol{\lambda}$, the conditional distribution of $\boldsymbol{\lambda}' \boldsymbol{\xi}_n$ given $\mathbf{X}_n$ is a CD for $\boldsymbol{\lambda}' \boldsymbol{\theta}$ where the $k \times 1$ random vector $\boldsymbol{\xi}_n$ has the distribution $H_n(\cdot)$ given $\mathbf{X}_n$.

Using the definition of asymptotic CD on the real line, one has a natural extension of the above definition to the asymptotic version. With this definition, for example, the raw bootstrap distribution in $\mathbb{R}^k$ remains an asymptotic CD, under asymptotic symmetry.

An $l$-CD $H_{1,n}$ is more precise than $H_{2,n}$ (both for the same $\boldsymbol{\theta}$ in $\mathbb{R}^k$), if the CD for $\boldsymbol{\lambda}'\boldsymbol{\theta}$ given by $H_{1,n}$ is more precise than that given by $H_{2,n}$, for all $k \times 1$ vectors $\boldsymbol{\lambda}$. For the normal mean vector $\boldsymbol{\mu}$, with known dispersion $\Sigma$, $1 - \Phi\big(\sqrt{n}\Sigma^{-\frac{1}{2}}(\bar{\mathbf{X}}_n - \boldsymbol{\mu})\big)$ is the most precise CD for $\boldsymbol{\mu}$.

Let $A_n(\hat{\boldsymbol{\theta}} - \boldsymbol{\theta})$ have a completely known absolutely continuous distribution $G(\cdot)$, where $A_n$ is non-singular, non-random matrix. Then $H_n(\cdot)$ defined by $H_n(\boldsymbol{\theta}) = 1 - G\big(A_n(\hat{\boldsymbol{\theta}} - \boldsymbol{\theta})\big)$ is an $l$-CD for $\boldsymbol{\theta}$. If $G$ is the limiting distribution then the above $H_n$ is an asymptotic CD, in which case $A_n$ can be data-dependent. A useful property of this extension to $\mathbb{R}^k$ is the fact that the CD content of a region behaves like $p$-values (in limit), as it does in the real line case. See, Theorem 4.2 of Liu and Singh [18] for the special case of bootstrap distribution in this context.

It is evident from the definition that one can obtain a CD for any linear function $A\boldsymbol{\theta}$ from that of $\boldsymbol{\theta}$, by obtaining the distribution of $A\boldsymbol{\xi}_n$. The definition also entails that for a vector of linear functions, the derived joint distribution (from that of $\boldsymbol{\xi}_n$) is an $l$-CD. For a non-linear function though, one can in general get an asymptotic CD only. Let $H_n(\cdot)$ be an asymptotic $l$-CD for $\boldsymbol{\theta}$. Suppose random vector $\boldsymbol{\xi}_n$ follows the distribution $H_n$. Consider a possibly non-linear function $g(\boldsymbol{\theta}) : \mathbb{R}^k \to \mathbb{R}^\ell, \ell \leq k$. Let each coordinate of $g(\boldsymbol{\theta})$ have continuous partial derivatives in a neighborhood of $\boldsymbol{\theta}_0$. Furthermore, suppose the vector of partial derivative at $\boldsymbol{\theta} = \boldsymbol{\theta}_0$ is non-zero, for each coordinate of $g(\boldsymbol{\theta})$.

**Theorem 5.1.** *Under the setup assumed above, the distribution of $g(\boldsymbol{\xi}_n)$ is an asymptotic $l$-CD of $g(\boldsymbol{\theta})$, at $\boldsymbol{\theta} = \boldsymbol{\theta}_0$, provided the $H_n$ probability of $\{||\boldsymbol{\theta} - \boldsymbol{\theta}_0|| > \epsilon\}$ $\xrightarrow{P} 0$, for any $\epsilon > 0$.*

*Proof.* The results follows from the following Taylor's expansion over the set $||\boldsymbol{\theta} - \boldsymbol{\theta}_0|| \leq \epsilon$: $g(\boldsymbol{\xi}_n) = g(\boldsymbol{\theta}_0) + \Delta(\boldsymbol{\xi}_n)(\boldsymbol{\xi}_n - \boldsymbol{\theta}_0)$, where $\Delta(\boldsymbol{\xi}_n)$ is the matrix of partial derivative of $g(\cdot)$ at $\boldsymbol{\theta}$ lying within the $\epsilon$-neighborhood of $\boldsymbol{\theta}_0$. □

**Remark 5.1.** Given a joint $l$-CD for $\boldsymbol{\theta}$, the proposition prescribes an asymptotic method for finding a joint $l$-CD for a vector of functions of $\boldsymbol{\theta}$, or for just a single function of $\boldsymbol{\theta}$. This method can inherit any skewness (if it exists) in the $l$-CD of $g(\boldsymbol{\theta})$. This will be missed if direct asymptotics is done on $g(\hat{\boldsymbol{\theta}}) - g(\boldsymbol{\theta}_0)$.

On the topic of combining joint $l$-CDs, one natural approach is by using the univariate CDs of linear combination, where combination is carried out by the methods discussed in Singh, Xie and Strawderman [26] . The problem of finding a combined joint $l$-CD comes down to finding a joint distribution agreeing with the univariate distributions of linear combinations, if it exists. The existence problem (when it is not obvious) could perhaps be tackled via characteristic functions and Bochner's Theorem. It may also be noted that an asymptotic combined $l$-CD of a nonlinear function of the parameters can be constructed via Theorem 5.1 and the methods of Singh, Xie and Strawderman [26].

### *5.2. Confidence distribution and data depth*

Another requirement on a probability distribution $H_n$ on $\mathbb{R}^k$ (based on a data set $\mathbf{X}_n$) to be a confidence distribution for $\boldsymbol{\theta}$ (a $k$-column vector) should naturally be: the $100t\%$ "central region" of $H_n$ are confidence regions for $\boldsymbol{\theta}$, closed and bounded, having the coverage level $100t\%$. We define such a $H_n$ to be a c-CD, where $c$ stands for circular or central. We note in Remark 5.2 that the notions of $l$-CD and c-CD match in a special setting.

**Definition 5.2.** A function $H_n(\cdot) = H_n(\cdot, \mathbf{X}_n)$ on $\Theta \in \mathbb{R}^k$ is called a Confidence Distribution in the circular sense (c-CD) for $k \times 1$ multiparameter $\boldsymbol{\theta}$, if (i) it is a probability distribution function on the parameter space $\Theta$ for each fixed sample set $\mathbf{X}_n$, and (ii) the $100t\%$ "central region" of $H_n(\cdot)$ is a confidence region for $\boldsymbol{\theta}$, having the coverage level $100t\%$ for each $t$.

By central regions of a distribution $G$, statisticians usually mean the elliptical regions of the type $y' \sum_G^{-1} y \leq a$. This notion of central regions, turns out to be a special case of central regions derived from the notion of data-depth. See Liu, Parelius and Singh [16] among others, for various concepts of data-depth and more. The elliptical regions arise out of so-called Mahalanobis-depth (or distance). The phrase data-depth was coined by J. Tukey. See Tukey [27], for the original notion of Tukey's depth or half-space depth. In recent years, the first author, together with R. Liu, has been involved in developing data-depth, especially its application. For the reader's convenience, we provide here the definition of Tukey's depth. $TD_G(\mathbf{x}) = \min P_G(H)$, where the minimum is over all half-spaces $H$ containing $\mathbf{x}$. On the real line, $TD_G(\mathbf{x}) = \min \big( G(x), 1 - G(x) \big)$. A notion of data-depth is called affine-invariant if $D(\mathbf{X}, \mathbf{x}) = D(A\mathbf{X} + \mathbf{b}, A\mathbf{x} + \mathbf{b})$ where $D(\mathbf{X}, \mathbf{x})$ is the depth of the point $\mathbf{x}$ w.r.t. the distribution of a random vector $\mathbf{X}$. Here $A\mathbf{X} + \mathbf{b}$ is a linear transform of $\mathbf{X}$. The above mentioned depths are affine-invariant. For an elliptical population, the depth contours agree with the density contours (see Liu and Singh [17]).

For a given depth $D$ on a distribution $G$ on $\mathbb{R}^k$, let us define the centrality function
$$C(\mathbf{x}) = C(G, D, \mathbf{x}) = P_G\{\mathbf{y} : D_G(\mathbf{y}) \leq D_G(\mathbf{x})\}.$$

Thus, the requirement (ii) stated earlier on $H_n$ in Definition 5.2, can be restated as: $\{\mathbf{x} : C(H_n, D, \mathbf{x}) \geq \alpha\}$ is a $100(1-\alpha)\%$ confidence region for $\boldsymbol{\theta}$, for all $0 \leq \alpha \leq 1$. This is equivalent to

Requirement (ii)': $C(\boldsymbol{\theta}_0) = C(H_n, D, \boldsymbol{\theta}_0)$, as a function of sample set $\mathbf{X}_n$, has the $U[0,1]$ distribution.

For a c-CD $H_n$, let us call the function $C(\cdot) = C(H_n, D, \cdot)$ a confidence centrality function (CCF). Here $D$ stands for the associated depth.

Going back to the real line, if $H_n$ is a CD, one important CCF associated with $H_n$ is
$$C_n(x) = 2\min\{H_n(x), 1 - H_n(x)\}.$$

This CCF gives rise to the two-sided equal tail confidence intervals. The depth involved is the Tukey's depth on the real line.

Next, we present a general class of sample dependent multivariate distributions which meet requirement (ii)'. Let $A_n(\hat{\boldsymbol{\theta}} - \boldsymbol{\theta})$ have a cdf $G_n(\cdot)$ independent of parameters. The nonsingular matrix $A_n$ could involve the data $\mathbf{X}_n$. Typically $A_n$ is a square root of the inverse dispersion matrix of $\hat{\theta}$ (or an estimated dispersion matrix). Let $\eta_n$ be an (external) random vector, independent of $\mathbf{X}_n$, which has the cdf $G_n$. Let $H_n$ be the cdf of $\hat{\boldsymbol{\theta}} - A_n^{-1}\eta_n$ for a given vector $\mathbf{X}_n$.

**Theorem 5.2.** *Let $D$ be an affine-invariant data-depth such that the boundary sets $\{\mathbf{x} : C(G_n, D, \mathbf{x}) = t\}$ have zero probability under $G_n$ (recall the centrality function $C(\cdot)$). Then $H_n(\cdot)$ defined above meets Requirement (ii)' and it is a c-CD.*

*Proof.* For any $t \in (0,1)$, in view of the affine-invariance of $D$, the set

$$\{\mathbf{X}_n : C(H_n, D, \boldsymbol{\theta}_0) \leq t\} = \{\mathbf{X}_n : C(G_n, D, A_n(\hat{\boldsymbol{\theta}} - \boldsymbol{\theta}_0)) \leq t\}$$
$$= \{X_n \text{ corresponding to } A_n(\hat{\theta} - \theta_0) \text{ lying in the outermost } 100t\%$$
$$\text{of the population } G_n\}.$$

Its probability content $= t$, since $A_n(\hat{\theta} - \theta)$ is distributed as $G_n$. □

**Remark 5.2.** The discussion preceeding Theorem 5.1 and the result in Theorem 5.2 imply that this $H_n$ is both a $l$-CD and a c-CD of $\boldsymbol{\theta}$ when $A_n$ is independent of the data. The $l$-CD and c-CD coincide in this special case! When $A_n$ is data-dependent, one has an exact c-CD, but only an asymptotic $l$- CD.

Given two joint c-CDs $H_{1n}$ and $H_{2n}$, based on the same data set, their precision could be compared using stochastic comparison between their CCFs involving the same data-depth. More precisely, let $C_1, C_2$ be the CCFs of two joint CDs $H_{1n}$ and $H_{2n}$, induced by the same notion of data-depth, i.e., $C_i(\mathbf{x}) = C(H_{in}, D, \mathbf{x}) = $ {fraction of the $H_{in}$ population having $D$-depth less than or equal to that of $\mathbf{x}$}. One would define: $H_{1n}$ is more precise than $H_{2n}$ when $\theta = \theta_0$ prevails, if $C_1(\mathbf{x}) \overset{\text{sto}}{\leq} C_2(\mathbf{x})$, under $\boldsymbol{\theta} = \boldsymbol{\theta}_0$, for all $\mathbf{x} \neq \boldsymbol{\theta}_0$.